\documentclass[12pt,leqno,a4paper]{article}
\usepackage{hyperref}
\usepackage{amsmath}
\usepackage{amssymb}
\usepackage{amsthm}
\usepackage{graphicx,epstopdf}
\usepackage{epsfig,graphics}

\hypersetup{
    colorlinks=true,
    linkcolor=black,
    filecolor=black,      
    urlcolor=blue,
    citecolor=red,
}

\makeatletter
\newcommand{\guio}[1]{\nobreakdash-\hspace{0pt}}

\newtheorem*{theorem*}{Theorem}

\newtheorem*{lemmam*}{Lemma 5. (\cite{MOV})}
\newtheorem*{lemma*}{Lemma}
\newtheorem*{corollary*}{Corollary}

\newtheorem{remark*}{Remark}
\newtheorem*{conjecture}{Semiadditivity of $\alpha$}

\theoremstyle{definition}

\newtheorem*{acknowledgements}{Acknowledgements}

\newcommand{\Rn}{\mathbb{R}^n}
\newcommand{\R}{\mathbb{R}}
\newcommand{\C}{\mathbb{C}}

\newcommand{\pv}{\operatorname{p.v.}}

\newcommand{\ep}{\varepsilon}
\newcommand{\e}{\epsilon}

\newcommand{\bm}{\operatorname{BMO}}

\title{Birth and life of the $L^{2}$ boundedness of the Cauchy Integral on Lipschitz graphs}

\author{Joan \ Verdera}

\date{}

\begin{document}
\maketitle
\begin{abstract}
We review various motives for considering
the problem of estimating the Cauchy Singular Integral on Lipschitz
graphs in the $L^{2}$ norm. We follow the thread that led to the solution and then describe a few of the innumerable applications
and ramifications of this fundamental result. We concentrate on its influence in complex analysis, geometric measure theory
and harmonic measure.

\bigskip


\medskip

\end{abstract}

\section{Introduction}
In 1982 Coifman, McIntosh and Meyer proved that the Cauchy Singular Integral on a Lipschitz graph is $L^2$ bounded with respect
to arc length on the curve \cite{CMM}. This is a deep result, simple to state, elegant, direct. In spite of its apparently specialized nature, 
it lies at the core of subtle important problems in PDE, and complex and 
harmonic analysis. It had resisted several attempts, devised by excellent mathematicians,  to discover a way to the proof. 
The story begins with Calder\'on and his first commutator, which he proved to be $L^2$ bounded in 1965 \cite{C1}. It continues when an important
breach was opened in 1977 
by Calder\'on himself, who showed the estimate assuming that the Lipschitz graph has uniformly small slopes everywhere \cite{C2}. This triggered a general mobilisation for the definitive assault, successfully achieved by Coifman, McIntosh and Meyer.

The impact of the CMcM (Coifman-McIntosh-Meyer) Theorem has been impressivily gigantic in a variety of areas, like PDE in Lipschitz domains, Kato square root 
problem, harmonic measure, complex analysis, singular integrals in non-homogeneous spaces and many more. The plan of this paper is, 
in the first place, to describe how the Russian School (Privalov, Vitushkin) walked the path to the walls of the CMcM Theorem. Their
motivation was very different from Calder\'on's, which has been explained very clearly in various excellent expositions \cite{M1,M2, S}, and closer
to my personal mathematical interests. Then I will describe some of the fantastic realms the magic key provided by CMcM opens the way:
Painlev\'e's problem, semi-additivity of analytic capacity, rectifiability, harmonic measure and many others.
The exposition does not aim at completeness, which is practically impossible in a short paper, and is very personal in all of its
aspects.

Let $A : \R \rightarrow \R$ be a Lispchitz function with Lipschitz constant 
$$\|A\|_{\text{Lip}}= \sup_{x \neq y} \frac{|A(y)-A(x)|}{|y-x|}$$
and let $\Gamma = \{(x,A(x)): x \in \R\}$ be its graph. The Cauchy Singular Integral is
\begin{equation}\label{csi}
 C(f)(z)= \pv \int_\Gamma \frac{1}{z-w}\, f(w)\,dw, \quad z=x+i A(x),
\end{equation}
where the principal value 
\begin{equation*}\label{}
\pv \int_\Gamma \frac{1}{z-w} \,f(w)\,dw = \lim_{\e \to 0}\int_{w\in \Gamma : |w-z|>\e} \frac{1}{z-w} \,f(w)\,dw.
\end{equation*}
 exists \text{a.e.} with respect to arc length on the graph for $f$ smooth with compact support, by an elementary argument.

Then the ``key that opens all doors'',  in the words of Y.Meyer, is the following inequality

\begin{theorem*}\label{cau}
There exists a constant $C$ depending only on $\|A\|_{\operatorname{Lip}}$ such that
\begin{equation}\label{caul2}
 \|C(f)\|_{L^2(\Gamma)} \le C\, \|f\|_{L^2(\Gamma)}, \quad f \in L^2(\Gamma).
 \end{equation}
\end{theorem*}

We are implicitly stating that one first shows the inequality for $f\in C^{\infty}_0(\Gamma)$ and then applies the standard machinery to show that the principal value integral in \eqref{csi} exists \text{a.e.} for each function in $L^2(\Gamma)$ and the $L^2$ estimate \eqref{caul2}   holds.

There is another equivalent expression for the Cauchy Singular Integral. Since for $w(y)=y+iA(y)$ one has $dw = (1+i A'(y))\,dy,$ 
incorporating the bounded factor $1+i A'(y)$ into $f(y+i A(y)),$ on can consider the equivalent operator

\begin{equation*}\label{}
\begin{split}
C(f)(x) &= \pv \int_\Gamma \frac{1}{x-y+i(A(x)-A(y))} \,g(y)\,dy \\*[8pt]
& = \lim_{\e \to 0}\int_{|y-x|>\e} \frac{1}{x-y+i (A(x)-A(y))} \,g(y)\,dy,
\end{split}
\end{equation*}
where now $g(y) = f(y+iA(y)) (1+iA'(y))$ with $f \in C_0^\infty(\Gamma).$  The set of $g$ is dense in $L^2(\R)$ and the inequality to be proven is
\begin{equation}\label{caul2r}
 \|C(g)\|_{L^2(\R)} \le C\, \|g\|_{L^2(\R)}.
\end{equation}
In the first formulation the ambient space is the Lipschitz graph and the measure arc length on the graph, which is clearly doubling
(the measure of a ball centred at $a$ of radius $2r$ is less than a 
constant times the measure of the ball with the same center and radius $r$). In the second formulation the ambient space is $\R$ and 
the underlying measure the familiar one dimensional Lebesgue measure. The main obvious difficulty for \eqref{caul2r}
is that the kernel is not of convolution type, except for $A$ linear (in which case you get a multiple of the Hilbert transform), and therefore the Fourier transform seems 
to be a forbidden tool. In spite of this there
is a proof of \eqref{caul2r} that relies on a reduction to a Fourier transform computation \cite{MV}.

The problem is truly beautiful and challenging. One feels immediately (and naively) that one can try something. I remember myself
considering polygonal graphs and trying to make estimates depending only on the biggest slopes. At another time I tried to emulate 
Loomis \cite{L} and I found myself making calculations with linear combinations of point masses. In both occasions I did not get anywhere.

\section{The fifties: Privalov}

Let $\Gamma$ be a closed rectifiable Jordan curve enclosing a bounded domain $D.$ Privalov was interested in non-tangential limits of
the Cauchy integral of a integrable function with respect to arc length on $\Gamma.$ The Cauchy Integral of such an $f$ is defined on 
$\C \setminus \Gamma$ by
\begin{equation*}\label{}
 C(f)(z) = \frac{1}{2\pi i}\int_\Gamma \frac{1}{w-z}\, f(w)\,dw, \quad z \notin \Gamma.
\end{equation*}
One would like to take 
non-tangential limits of $C(f)(z)$ at points $a \in \Gamma$ where a tangent exists, that is, limits when the variable $z$ approaches
the point $a$ either from $D$ or from the exterior, but is constrained to remain in a double sided cone $K(a)$ with axis the normal line 
to $\Gamma$ at $a$ and some fixed aperture.
One writes
\begin{equation*}\label{}
 C(f)^+(a) = \lim_{ K(a)\cap D \ni z \to a} C(f)(z), 
\end{equation*}
and
\begin{equation*}\label{}
 C(f)^-(a) = \lim_{ K(a)\cap \overline{D}^c \ni z \to a} C(f)(z),
\end{equation*}
whenever these limits exist. There is a third operator involved, the Singular Cauchy Integral, defined in terms of principal values,
provided they exist, namely
\begin{equation*}\label{}
 C(f)(a) = \frac{1}{2\pi i} \pv \int_\Gamma \frac{1}{w-a} \,f(w)\,dw, \quad a \in \partial D.
\end{equation*}
Privalov was interested in deciding whether the limits above are defined \text{a.e.}\ with respect to arc length on $\Gamma.$ He proved in
\cite{P} that, given $f \in L^1(\Gamma),$ at arc length almost all points of $\Gamma$ the existence of one of these limits implies the existence of the other two
 and one has the classical $Plemelj-Sohotski$ formulas 
\begin{equation*}\label{}
\begin{split}
 C(f)^+(a) & = C(f)(a) + \frac{1}{2} f(a), \\*[8pt]
  C(f)^-(a) & = C(f)(a) - \frac{1}{2} f(a), \\*[8pt]
  C(f)^+(a) &- C(f)^-(a)  = f(a).
\end{split}
\end{equation*}
\vspace{0.2cm}

Privalov did not mention in his book any connection with $L^2$ estimates for the Cauchy Singular Integral, in spite of the fact that
the reduction of the \text{a.e.}\ existence of these limits to Calder\'on's Theorem on $L^2$ estimates on Lipschitz graphs with small 
constant is rather easy. Havin showed in 1965 \cite{H} that for the \text{a.e.} existence problem one can reduce the case of integrable functions to that of continuous 
functions.
Much later Dynkin modified slightly Havin's argument to make the reduction from a general rectifiable curve to a Lipschitz curve with small constant.
The interested reader will find more details in Dynkin's excellent survey in the Encyclopedia of Mathematics \cite[p. 216]{Dy}.

Therefore one is led to questions about $L^2$ boundedness by issues concerning \text{a.e.}\ existence of certain limits. It is somehow
odd that the problem was not pushed by the Russian school to the very heart of the matter, namely, inequality \eqref{caul2}.

\section{The sixties: Vitushkin}

In the sixties the problem of uniform rational approximation was very popular in the United States, Canada and the Soviet Union. Given a compact 
set $K$ and a continuous function $f$ on $K$ one would like to know under which conditions $f$ can be uniformly approximated on $K$
by rational functions with poles off $K,$ or, equivalently (by Runge's Theorem), by functions analytic on neighbourhoods of $K.$ Assume 
that $f$ has been extended to a continuous function on $\C.$ Then Vitushkin proved in a remarkable paper \cite{Vi} that a necessary and sufficient
condition for the approximability of $f$ by rational functions without poles on $K$ is that there exists a function 
$\omega(\delta), \; \delta >0,$ tending to $0$ with $\delta,$ such that 
for each open square $Q$ of side length $\delta$ one has
\begin{equation}\label{Vi}
\left| \int_{\partial Q} f(z) \, dz \right| \le \omega(\delta)\,  \gamma(Q\setminus K).
\end{equation}
Here $\gamma$ is a set function, called analytic capacity, introduced by Ahlfors in \cite{A}. To understand better the above inequality
it is necessary to devote some lines to $\gamma.$ The analytic capacity of a compact set $E$ is
\begin{equation}\label{ancap}
\gamma(E)= \sup |h'(\infty)|,
\end{equation}
where the supremum is taken over all bounded analytic functions $h$ on $\C \setminus E$ satisfying the normalization condition
$|h(z)|\le 1, \; z \notin E.$ If $F$ is any set then $\gamma(F)$ is the supremum of $\gamma(E)$  on compact sets $E\subset F.$
It is straightforward that $\gamma(E)=0$ if and only if $E$ is removable for bounded analytic functions, so that
$\gamma$  quantifies the notion of non-removability for $H^\infty.$   For example, the analytic capacity of a disc $B(a,r)$ 
of center $a$
and radius $r$ is exactly the radius. Painlev\'e proved
that a set of zero Hausdorff length is removable and it is not difficult to argue that a set of Hausdorff dimension greater than
$1$ is non-removable. There are inequalities in terms of $\gamma$ that quantify the above statements, namely
\begin{equation*}\label{}
 C_\beta \, H^{\beta}_\infty(E)^{\frac{1}{\beta}} \le   \gamma(E) \le  H^1_\infty(E),\quad 1 < \beta,
\end{equation*}
where $H^\beta_\infty$ stands for $\beta$-dimensional Hausdorff content.
Ahlfors asked in \cite{A} for geometric characterizations of removability, and this has been called
since then the Painlev\'e problem. At that time papers were written without a list of references, so that no paper by 
Painlev\'e was mentioned by Ahlfors. I believe that  
Ahlfors himself was the first to ask for geometric characterizations of $H^\infty$ removability. He wrote ``I have not been able to push the solution so far''. 
Tolsa found the solution in \cite{T2}, 54 years after Ahlfors' paper, building on hard brilliant work by the previous generations.
In particular, the proof is rooted in the ground of the magic key \eqref{cau}. 

Let us discuss Vitushkin's condition \eqref{Vi}. If $Q\subset K$ then the right hand side of \eqref{Vi} vanishes and $f$ is analytic on the interior
of $K$ by Morera's theorem. If $Q$ does not intersect $K$ then \eqref{Vi} is  trivially satisfied with $\omega$ the modulus
of continuity of $f$ (substract from $f$ its value at the center of $Q$). Thus the condition really involves squares that intersect $\partial K$ and should be viewed as a weak analyticity
condition for $f$ on the boundary of $K.$ It can be easily seen that \eqref{Vi} implies Mergelyan's theorem : if the complement of $K$
is connected then any continuous function on $K,$ analytic on its interior, can be uniformly approximated on $K$ by holomorphic polynomials. 

In fact, there is a close relative of $\gamma,$ called continuous analytic capacity, that is more tightly connected to approximation
issues than $\gamma.$ The continuous analytic capacity $\alpha(E)$ of a compact set $E$ is the supremum in \eqref{ancap} taken over
analytic 
functions on $\C \setminus E,$ that extend continuously to $\C$ and satisfy the normalization $|h(z)|\le 1, \; z \notin E.$
If $F\subset\C,$ then $\alpha(F)$ is the supremum of $\alpha(E)$ over compact subsets $E$ of $F.$ If the function $f$ is continuous
on the plane and analytic on the interior of $K,$ then one has a weaker form of \eqref{Vi}, namely,
\begin{equation}\label{Vialfa}
\left| \int_{\partial Q} f(z) \, dz \right| \le \omega(\delta)\,  \alpha(Q\setminus \mathring{K}),
\end{equation}
where $\omega$ is the modulus of continuity of $f.$  It is easily seen that $\gamma $ and $\alpha$ agree on open sets, so that in right hand side of
\eqref{Vi} $\gamma (Q\setminus K)$ can be replaced by $\alpha (Q\setminus K).$  Then a consequence of \eqref{Vi} is  
the following striking result \cite{Vi}.

\begin{theorem*}\label{vith}
 For a compact subset $K$ of the plane the following are equivalent.

(i) Every continuous function on $K,$ analytic on its interior,
can be uniformly approximated on $K$ by rational functions with poles off $K.$

(ii) For each open square $Q$
\begin{equation}\label{alfac}
      \alpha(Q\setminus \mathring{K}) \le C \, \alpha(Q\setminus K),
\end{equation}
where the constant $C$ is independent of $Q.$
\end{theorem*}

In Vitushkin's terminology a singularity of a function $f$ is a point such that in no neighbourhood of the point $f$ is 
analytic. In other words, it is a point in the support of the distribution $\overline{\partial} f.$
The proof that condition \eqref{alfac} implies the approximation property on $K$ consists in localizing the singularities of the function
to be approximated and then pushing them off $K$ by \eqref{alfac}. Condition \eqref{alfac} is explicit enough to provide several
examples of sets $K$ with the approximation property. For example, if the complement of $K$ has finitely many components then
\eqref{alfac} holds. It also holds if the inner boundary of $K$ is very small. The inner boundary of $K$ is the set of boundary points which
are not in the closure of a bounded component of the complement of $K.$ Imagine a compact set constructed from the closed unit disc
by deleting a sequence of mutually disjoint open discs accumulating at the origin. Then the origin is the only point in the inner
boundary and property \eqref{alfac} holds. Vitushkin proved that if these discs accumulate to a Jordan arc of class 
$C^{1+\varepsilon}$ then \eqref{alfac} still holds, in spite of the fact that the inner boundary is a one
dimensional object.
The drawback was that his method did not cover the case $\varepsilon=0.$ The situation was then, let us say, uncomfortable.

\begin{figure}[h]
\begin{center}
\includegraphics{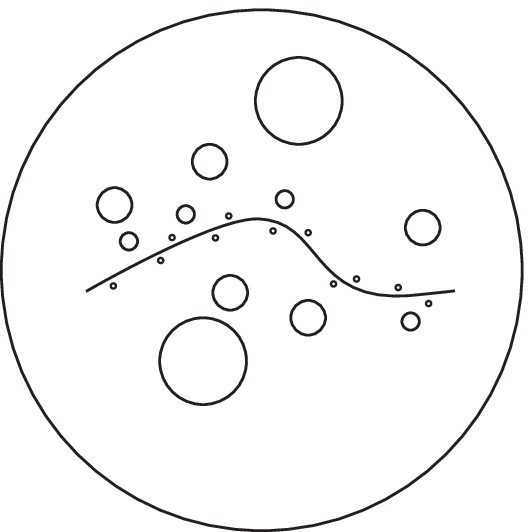}
\end{center}
\end{figure}

If $J$ denotes the arc, $E$ the union of the boundaries of the deleted discs and $Q$ is an open square, then 
$$Q\setminus \mathring{K} = (Q \cap J)\cup (Q\cap E) \cup (Q\setminus K)$$
and Vitushkin's approximation technique gives readily that
$$
\alpha(Q\setminus \mathring{K}) = \alpha ((Q \cap J) \cup (Q\setminus K)).
$$
It is not difficult to show that a set of finite length (one dimensional Hausdorff measure) has zero continuous analytic capacity.
Thus the following conjecture arises.
\begin{conjecture}
There exist a constant $C$ such that
$$\alpha(K_1\cup K_2) \le C \left(\alpha(K_1)+\alpha(K_2)\right),$$
for all compact sets $K_1$ and $K_2.$
\end{conjecture}
This, of course, implies that  the approximation property holds if the inner boundary is a  $C^1$ arc and even if it is 
a rectifiable arc. There is an analog of the semiadditivity conjecture for analytic capacity.

It turns out that both forms of semiaddivity were proved to hold true many years later by Tolsa in \cite{T2} and \cite{T3} using sophisticaded
methods based essentially on a special form of the $T(b)$ theorem, which is a close descendant of the CMcM theorem. Again, as we said
when referring to the solution of Painlev\'e's problem, Tolsa's results depend on previous subtle, outstanding work by other mathematicians
(David, Jones, Nazarov, Melnikov, Mattila, Treil, Volberg and others).

One fact worth mentioning is that the smoothness barrier Vitushkin found when applying his methods was $C^1$ 
(or the Lipschitz condition for that matter). 
The analogy with Calder\'on's Lipschitz graph problem is not only formal. It is the same kind of obstruction. The context is not so
transparent as with the graph, but it was rather clear in the eighties that both issues were intimately related \cite{V1}.

To close this section I would like to mention a wonderful result of Mazalov \cite{Mz}. He proved that for any compact
$K$ each continuous function  $f$ on $K$ satisfying the equation $\overline{\partial}^2 f=0$ on $\mathring{K}$ can be uniformly 
approximated on $K$ by functions of the type $r_0(z) + \overline{z}\, r_1(z),$ where $r_0$ and $r_1$ are rational functions with poles
off $K.$ Thus there is no capacitary restriction on the approximation property for the operator $\overline{\partial}^2,$ which turns out to be universal. This is due to
the fact that the fundamental solution of $\overline{\partial}^2,$ which is $\overline{z}/\pi z,$  is bounded. The main ingredient of the proof, besides a deep elaboration of Vitushkin's techniques,
is an intricate clever construction of a Lipschitz graph which allows the use of the magic key \eqref{cau}.

\section{Big pieces of Lipschitz graphs and the Denjoy conjecture}
In the mid eighties I went to Madrid for two consecutive years to listen to Yves Meyer, who had been invited there to lecture. 
I remember these lectures very vividly, because he was talking about results which seemed to me of the utmost importance for what 
I was trying to do. Also because the speaker was superb, always clear, a soft, contained, but perceptible passion permeating the exposition.

I learnt there that a student of his, by the name of Guy David, had shown that  rectifiable curves $\Gamma$ for which the Cauchy
Integral is $L^2(\Lambda)$ bounded, $\Lambda$ being arc length on the curve, were those satisfying the Ahlfors
condition
\begin{equation*}\label{ahl}
 \Lambda(B(z,r)) \le C_0 \, r, \quad z \in \Gamma, \quad 0 < r < \operatorname{diam} (\Gamma),
\end{equation*}
where $C_0$ a positive constant independent of $z$ and $r.$  The inequality from below 
$C_0^{-1} \, r \le \Lambda(B(z,r))$
is automatic, because $\Gamma$ is a curve. Thus the Ahlfors condition is indeed equivalent to
\begin{equation}\label{ahldavid}
C_0^{-1} \, r \le \Lambda(B(z,r)) \le C_0 \, r, \quad z \in \Gamma, \quad 0 < r < \operatorname{diam} (\Gamma),
\end{equation}
which is now called the AD (Ahlfors-David) condition.
The proof of David's theorem, extremely elegant, 
consists in finding a ``big piece'' of a Lipschitz graph inside each ball centered at a point on the curve. More precisely, for each
ball
$B=B(a,r), \; a \in \Gamma,$ there exists a Lipschitz graph $G=G_B$, possibly rotated, with Lipschitz constant controlled by $C$
such that $\Lambda(B(a,r)\cap \Gamma) \ge C^{-1}\, \Lambda(B(a,r)),$  where $C$ depends only on $C_0.$ The graph is constructed by applying the 
rising sun Lemma. On $G_B$ the Cauchy Integral is $L^2$
bounded with respect to $\Lambda$ with a constant depending only on the Ahlfors constant.  This situation is ideal for a good $\lambda$ 
inequality, which completes the argument \cite{D1}. In fact, the local Lipschitz graph can be taken so that $\|A'\|_\infty$ is as small
as desired by a recursive reasoning. Thus the $1982$ CMcM Theorem can be reduced to Calder\'on's 1977 small Lipschitz constant theorem.

What interested me a lot at that time is that David's theorem implies a quantitative version of the Denjoy conjecture. 
Denjoy believed that he could show that a compact subset of a rectifiable curve  of positive length is non-removable for $H^\infty.$
Actually his argument worked only for the line, whence the conjecture.

In terms of analytic capacity, the Denjoy conjecture states that a compact set of positive length in a rectifiable curve has positive
analytic capacity.  We will see now that, given a curve
satisfying the AD condition \eqref{ahldavid},   there exists a constant $c$ depending 
only on the AD constant $C_0$ such that
\begin{equation}\label{denj}
 \gamma(K) \ge c \, \Lambda(K), \quad K \; \text{compact}\; \subset \Gamma.
\end{equation}
This is a quantitative solution of the Denjoy conjecture, which depends on the $L^2$ estimate
for the Cauchy Integral. This fact, together with a reduction argument to Lipschitz graphs with small constant, explains why the Denjoy
conjecture was a consequence of Calder\'on's 1977 result. Apparently Calder\'on was not aware at that time that he had solved the
Denjoy conjecture
and Donald Marshall wrote a short note (I remember a footnote : ``not intended for publication'') to describe the argument via 
the Garabedian $H^2$ description of analytic capacity \cite{Ga}. Later on a much more direct and suggestive way of obtaining
\eqref{denj} emerged, which I would like to discuss now.

Since the length measure on an Ahlfors curve satisfies the doubling condition (the measure of the double concentric ball is not greater
than a constant times the measure of the ball) standard CZ(Calder\'on-Zygmund) theory goes through to show the weak-type $(1,1)$ 
inequality for finite measures $\mu$ supported on $\Gamma$:
\begin{equation}\label{wt}
\Lambda \{z\in \Gamma : |C(\mu)(z)| > t \} \le C\, \|\mu\|, \quad t>0,
\end{equation}
where
\begin{equation*}
 C(\mu)(z)= \pv \int_\Gamma \frac{1}{z-w}\,d\mu(w), \quad z \in \Gamma.
 \end{equation*}
 As in classical CZ theory, one can get an $L^p$ inequality, $1<p<2,$ \,from \eqref{wt} by interpolating with the $L^2$ estimate; then
 by duality one reaches the full range $1<p<\infty.$ One even manages to have an $L^\infty-BMO$ estimate. But for \eqref{denj} one
 needs to end up with functions with bounded Cauchy Integral. 
It turns out that this can be achieved by dualizing the weak-type inequality \eqref{wt}, which sounds surprising at first glance, because
apparently nobody knows the dual of $L^{1,\infty}.$  That a weak-type inequality can be dualized was first discovered by  N.X. Uy,
a former student of Garnett, who proved that a compact 
set of positive area
in the plane is non-removable for Lipschitz holomorphic functions \cite{U}. Later Davie and Oksendal found a slightly
stronger, more convenient way of expressing the dual statement \cite{DO}. What they proved in half a page, via an ingenious
duality argument based on Hahn-Banach, is that the weak-type inequality  yields the following statement, which is indeed
equivalent to \eqref{wt}.
\begin{lemma*}
Given a compact subset $K$ of $\Gamma,$ there exist a measurable function $h$ supported on $K,$ $0\le h\le1,$ 
such that
$\Lambda(K) \le 2 \int h \,d\Lambda$ and 
\begin{equation*}\label{ib}
 |C(h\,d\Lambda)(z)| \le C, \quad \text{for each}\;z \notin K,
\end{equation*}
and for $\Lambda$ almost all $z \in \Gamma.$  The constant $C$ is comparable to that of \eqref{wt}.
\end{lemma*}
The function $C(h\,d\Lambda)$ is holomorphic and bounded off $K$ and the derivative at $\infty$ is $\int h \, d\Lambda,$
which gives readily \eqref{denj}. Therefore a non constant bounded holomorphic function off $K$ has been identified (provided
$K$ has non zero length).
Peter Jones has asked about a constructive argument for the existence of such a function. I heard once a statement about
the fact the the Hahn-Banach theorem for separable Banach spaces is already constructive. But anyway I understand Peter Jones' quest.

How did Uy prove the existence of a non linear holomorphic Lipschitz function off a compact set of positive area ? He dualized 
the weak inequality with respect to planar Lebesgue measure for the Beurling transform $B.$ He obtained a non zero bounded measurable 
function
$h$ supported on the compact set $K$ such that $B(h)$ is also bounded. Then he set $f=C(h \,dA)$ and he was done, because 
$\overline{\partial}f =h$ and $\partial f = B(h)$ are both bounded and hence $f$ is a Lipschitz function on the plane, analytic off $K$ but not on $K.$
It is the same argument that leads to \eqref{denj}. I took advantage
of Uy's theorem to settle the problem of $C^1$ approximation by rational functions in \cite{V0}.

It was clear to me at that time that the study of the Cauchy Integral as a CZ singular operator would play a decisive role in the 
understanding of analytic capacity, but I could not imagine  to what extent. Further progress into this direction was provided by
Yves Meyer in his lectures in Madrid the second year.

\section{David, Journ\'e and Semmes : T(1) and T(b)}

The success with the Cauchy Integral led naturally to consider the following question. Let $K(x,y)$ be a kernel defined off the diagonal of
$\R^n$
satisfying the usual growth and smoothness conditions of CZ theory:
\begin{equation*}\label{}
\begin{split}
 |K(x,y)| &  \le \frac{C}{|x-y|^n}, \quad x \neq y,\\*[5pt]
 |\nabla K(x,y)| & \le \frac{C}{|x-y|^{n+1}}, \quad x \neq y.
\end{split}
\end{equation*}
Let us restrict attention, to simplify matters, to antisymmetric kernels : $K(y,x)=-K(x,y).$ Consider the truncations
\begin{equation*}
 T_\e (f)(x) = \int_{|x-y|>\e} K(x,y) \,f(y)\,dy, \quad x \in \R^n, \quad \e >0,
\end{equation*}
which are well defined for $f$ in $L^p(\R^n), \; 1\le p< \infty$. We would like to find a necessary and sufficient condition so that the operator $T$ is bounded
on $L^2(\R^n)$, namely, so that
\begin{equation}\label{tl2b}
 \int |T_\e (f)(x)|^2\,dx \le C\,  \int |f(x)|^2\,dx,\quad f \in L^2(\R^n), \quad \e >0,
\end{equation}
with $C$ independent of $f$ and $\e.$ The preceding inequality for $f=\chi_B,$ \,$B$ a ball, is
\begin{equation}\label{balls}
 \int_B  |T_\e (\chi_B)(x)|^2\,dx \le C\,|B|, \quad \e >0,
\end{equation}
where $|E|$ stands for the Lebesgue measure of $E.$ The part of the integral off the ball satisfies the required inequality, as standard
reasoning shows. It is not immediate, but simple, to define $T_\e(1)$ (or $T(1)$, for that matter)
as a distribution modulo constants, that is, as a continuous functional on the space of compactly supported smooth functions
with zero integral. A classical argument derives from \eqref{balls} that $T_\ep(1)$ is in 
$\operatorname{BMO}(\R^n),$ with norm independent of $\ep.$  There is a small difficulty here related to the fact that the truncated 
kernel $\chi_{|x-y|>\e}(x,y) K(x,y)$ does not satisfy
the smoothness condition involving the gradient, which is overcome by remarking that it satisfies Hormander's condition, 
uniformly in $\e$. One can also conclude that $T(1) \in \operatorname{BMO}(\R^n)$ without appealing to truncations.  
Two former students of Meyer, David and Journ\'e, proved the converse. Guy David says that Meyer suggested the use of paraproducts to complete the proof. The $T(1)$-Theorem reads as follows.

\begin{theorem*}[\bf{David and Journ\'e}]
The operator $T$ is bounded on $L^2(\R^n)$ if and only if $T(1)$ is in $\bm(\R^n).$
\end{theorem*}

 This can also be formulated by saying that \eqref{tl2b} and \eqref{balls} 
are equivalent. The outcome is that to check that an operator with standard antisymmetric kernel is $L^2(\R^n)$ bounded one only has to check 
the action of the operator on the function $1.$ This is an extremely powerful fact, as shown by the following examples.

The first Calder\'on commutator is
 \begin{equation}\label{commu}
C_1f(x)=\left(AH\frac{d}{dx}-H\frac{d}{dx}A\right)(f)=\text{p.v. }\int^{\infty}_{-\infty}\frac{A(y)-A(x)}{(y-x)^2}f(y)\,dy,
\end{equation}
where $A$ is a Lipschitz function on the real line and $H$ the Hilbert transfom.
Calder\'on's PDE motivation for tackling the problem of the $L^2(\R)$ boundedness of the first commutator has been explained in several
nice expository articles \cite{C1}, \cite{M1} and \cite{M2}. His celebrated proof was based on complex analytic methods involving the Hardy 
space $H^1$ on the upper half plane.  Now the first commutator
applied to the function $1$ turns out to be, after an integration by parts,
\begin{equation*}\label{}
\begin{split}
 \text{p.v. }\int^{\infty}_{-\infty}\frac{A(y)-A(x)}{(y-x)^2} \,dy & = 
 \text{p.v. }\int^{\infty}_{-\infty}\frac{1}{(y-x)} A'(y)\,dy = -H(A')(x),
\end{split}
\end{equation*}
which is a $\bm(\R)$ function because the Hilbert transform maps $L^\infty(\R)$ into $\bm(\R).$ 
The $T(1)-$Theorem yields immediately the  $L^2(\R)$ boundedness of the first commutator. 

Commuting with multiplication by $A$ several times leads the $n$-th order commutator, whose kernel is
$\frac{(A(y-A(x))^n}{(y-x)^{n+1}}.$ Coifman and Meyer had devoted a lot of energy for proving the sharpest possible $L^2(\R)$ estimates
for the higher order commutators, employing Fourier analysis methods \cite{CM}. The hope was to prove the Lipschitz graph Theorem in 
full generality.
The action of the $n$-th commutator on the function $1$ is precisely the action on the $(n-1)$-th commutator on $A'.$ Thus 
successive applications of the $T(1)$-Theorem lead to $L^2(\R)$ estimates for the $n$-th commutator with constant of the form
$C^n \|A'\|_\infty^n,$ with $C$ a numerical constant. This gives the Lipschitz graph Theorem for small Lipschitz constant,
the famous Calder\'on's result of 1977. The reason is that one can expand
the Cauchy kernel on the graph in terms of those of higher order commutators
\begin{equation*}\label{}
\begin{split}
\frac{1}{x-y+i(A(x)-A(y))} = \sum_{n=0}^\infty (-i)^n \frac{(A(x)-A(y))^n}{(x-y)^{n+1}}.
\end{split}
\end{equation*}
Hence the series of the corresponding operators can be summed up if $C \|A'\|_\infty < 1.$

The main drawback with the $T(1)$-Theorem at that time was that it could not be applied to get a proof of 
the magic key \eqref{cau}, just because you could not compute $T(1).$ Instead by Cauchy's Theorem you see immediately that
\begin{equation*}\label{}
C(1+i A')(x) = \pv \int_\Gamma \frac{1}{x-y+i(A(x)-A(y))} \,(1+i A'(y))\,dy =0.
\end{equation*}
Thus the bounded function $b:=1+i A'$ is mapped into the  $\bm(\R)$ function ``0'' and it is far from being zero, because
$\operatorname{Re}b = 1.$

David, Journ\'e and Semmes proved an extension of the $T(1)$- Theorem in which the function $1$ is replaced by a bounded function
$b$ which satisfies a non-triviality condition called para-accretivity. This means, in the $n$-dimensional context described before, 
that
$$
\left|\frac{1}{|B|} \int_B b(x) \, dx \right| \ge c >0, \quad \text{for all balls}\quad B.
$$
The $T(b)$-Theorem is the following.
\begin{theorem*}[\bf{David, Journ\'e and Semmes}]
The operator $T$ is bounded on $L^2(\R^n)$ if and only if $T(b)$ is in $\bm(\R^n)$ for some  bounded para-accretive function $b.$ 
\end{theorem*}

Accretivity refers to the property $\operatorname{Re}b(x) \ge \delta >0, \; x \in \R^n,$ obviously stronger than para-accretivity.
The $T(b)$-Theorem applied to $b=1+iA'$ yields the Lipschitz graph Theorem in full generality. 

The $T(b)$-Theorem is stated and proved by David, Journ\'e and Semmes on spaces of homogeneous 
type and for kernels
not necessarily antisymmetric. The underlying euclidean space is replaced by a metric space and the Lebesgue measure by a Borel 
doubling measure. One example of space of homogenous type is an Ahlfors regular curve with length as underlying measure.
However, one cannot apply $T(b)$ to Ahlfors regular curves, because the para-accretivity condition fails for the most
natural candidate, namely, the unit 
tangent vector. This is due to the fact  that the curve may intersect itself. The big pieces idea remains a fundamental contribution. $T(b)$ can be still applied to chord-arc curves.

Let us say a few words about a more recent and simpler proof (\cite{CJS} and \cite{D1b}). There are two steps. 

First you assume $T(b)=0$ and borrow ideas from the most basic wavelet,
the one leading to the Haar basis. You construct a Haar basis adapted to $b$ and you check that the off diagonal entries of the matrix of 
$T$ in this basis decay rapidly to $0.$ This gives boundedness via Schur's Lemma. 

The second step consists in finding an explicit operator $P$ bounded on $L^2(\R^n)$, 
called paraproduct, with the property that $P(b)=\beta,$ where $\beta$ is a given function in $\bm(\R^n).$ The paraproduct $P$ has a
kernel satisfying the usual growth conditions of CZ Theory, but smoothness fails. However, one still has the relevant consequences
of smoothness that allow to get $L^2(\R^n)$ boundedness.
Let $P$ be the paraproduct
associated with $\beta =T(b).$ Since $(T-P)(b)=0$ one concludes from the first part of the proof that $T-P$ is bounded. 
On the other hand,
$P$ is bounded by contruction and so the proof is complete.

It will become clear later on that dealing with analytic capacity requires to get rid of the doubling condition and still have
a $T(b)$-Theorem. That this can be done is a great success story that I will sketch briefly, but not in the next section, which is 
devoted to Menger curvature, a device that connects the analysis of the Cauchy Integral to the geometry of the triangle.

\section{The first commutator controls the Cauchy Integral}

It is very thrilling to read Meyer's expositions
on the sequence of events that lead from $L^2$ boundedness of the first commutator to that of the Cauchy Integral on Lipschitz graphs
\cite{M1} and \cite{M2}. Let us review this briefly.  Calder\'on proved that the first commutator is bounded in $L^2(\R)$ in 
1965 using complex variable methods and the Lusin area function in \cite{C1}. The Cauchy Integral was out of reach at that time.
Coifman and Meyer started working in the late seventies on
 $L^2$ estimates for higher order commutators by real variables methods. They achieved outstanding results
but not sharp enough to reach the Cauchy Integral. In 1977 Calder\'on proved the estimate on graphs with
small Lipschitz constant and finally Coifman, McIntosh and Meyer got the full result in 1982.  Applications to PDE on domains with 
$C^1$ boundary started immediately after
the small constant theorem \cite{FJR} and after the full result the flow of applications in PDE grew exponentially 
\cite{Ke}.
In particular, the solution of Kato square
root problem in dimension $2$ \cite{KM} was achieved. 

It is  shown in \cite{V4} that, indeed, the  Cauchy Integral is controlled by the first commutator. In other words, the 1965 result
of Calder\'on plus the $H^1-BMO$ duality imply the 1982 $\text{CMcM}$ Theorem. This is surprising and, in a certain sense,
a manifestation of the irony that pervades life. The first proof of an important result is rarely the best, or the shortest.
People go on working on the theme, new ways of facing the problem and new connections arise, new paths to the summit are discovered.
All this activity enlivens the subject, clarifies the results, makes them more accessible.
I would like to 
explain how one can prove $L^2$ boundedness of the first commutator and the Cauchy Integral in a relatively short way.

Let us start with the first commutator. The kernel is
\[
K(x,y)=\frac{A(y)-A(x)}{(y-x)^2}.
\]
 The truncations of the first
commutator are
\[
C_{1,\epsilon}(f)(x)=\int_{|y-x|>\epsilon}K(x,y)f(y)\,dy,\quad
x\in\mathbb{R}, \quad \e >0.
\]
Given an interval~$I$, we have
\begin{equation*}
\begin{split}
\int_IC_{1,\epsilon}(\chi_I)^2(x)\, dx &=\int_I\int_{I_{\epsilon}(x)}\int_{I_{\epsilon}(x)}K(x,y)K(x,z)\,dy\,dz \,dx\\
&=\iiint_{S_{\epsilon}}K(x,y)K(x,z)\,dx\,dy\,dz+O(|I|),
\end{split}
\end{equation*}
where 
$$I_{\epsilon}(x)=\{t\in I:|t-x|>\epsilon\}$$ 
and 
\[
S_{\epsilon}=\{(x,y,z)\in I^3: |y-x|>\epsilon,\,
|z-x|>\epsilon,\,|z-y|>\epsilon\}.
\]
This can be shown readily by splitting the triple integral on $I\times I_\e(x)\times I_\e(x)$ as the sum of three triple integrals
over
$$
U_\e =\{(x_1, x_2, x_3) \in I^3: |x_1-x_2|\le \e,\,|x_3-x_1| > 2 \e, \,|x_3-x_2|> \e \},
$$
$$
V_\e =\{(x_1, x_2, x_3) \in I^3: |x_1-x_2|\le \e,\, \e < |x_3-x_1| \le 2\e, \,|x_3-x_2|> \e \}
$$
and $S_\e$ respectively.

Note that the integrand in the triple integral over $S_\e$  is
already symmetric in $y$ and $z$. If we interchange $x$ by $y$, and $x$ by
$z$, we get two new different expressions for
$\int_IC_{1,\epsilon}(\chi_I)^2$. Therefore
\begin{equation}\label{compos}
 \int_IC_{1,\epsilon}(\chi_I)^2(x)\,dx=\frac{1}{3}\iiint_{S_{\epsilon}}S(x,y,z)\,dx\,dy\,dz+O(|I|)
\end{equation}
where
\[
S(x,y,z)=K(x,y)K(x,z)+K(y,x)K(y,z)+K(z,y)K(z,x).
\]
A straightforward computation yields
\[
S(x,y,z)=\left(\frac{\dfrac{A(y)-A(x)}{y-x}-\dfrac{A(z)-A(x)}{z-x}}{z-y}\right)^2,
\]
which is a non-negative quantity. Therefore we have expressed the $L^2$ norm on the left hand side of \eqref{compos}
in terms of a non-negative kernel, which allows us to avoid the fine cancellation effects which lie at the heart of CZ theory.

A beautiful simple Fourier transform estimate gives the following lemma.

\begin{lemma*}\label{fou}
If $a$ is a locally integrable function with derivative~$a'\in L^2(\mathbb{R})$, then
\[
\int^{\infty}_{-\infty}\int^{\infty}_{-\infty}\int^{\infty}_{-\infty} \left|\frac{\dfrac{a(y)-a(x)}
{y-x}-\dfrac{a(z)-a(x)}{z-x}}{z-y}\right|^2\,dx\,dy\,dz=c_0\int^{\infty}_{-\infty}|a'(x)|^2\,dx,
\]
$c_0$ being a numerical constant.
\end{lemma*}

The next step is to localize the $L^2$ identity above to the
interval~$I$ that was given.
This is a kind of technique one uses often when working with BMO. Indeed, in the next
inequality the BMO norm of $A'$ shows up naturally.

Apply the lemma to the function~$a=\chi_I(A-P_I)$, where
$P_I(x)=A'_I(x-\alpha)+A(\alpha)$, $A'_I=\frac{1}{|I|}\int_I A'$, $I=(\alpha,\beta)$.
The result is
\[
\iiint_{I^3}  \left(\frac{\dfrac{A(y)-A(x)}{y-x}-\dfrac{A(z)-A(x)}{z-x}}{z-y}\right)^2   \,dx\,dy\,dz \le C\int_I
|A'-A'_I|^2\le C\|A'\|_{\infty}^2|I|,
\]
and thus
\begin{equation}\label{eq4.2bis}
\int_I|C_{1,\epsilon}(\chi_I)|^2\le C(1+\|A'\|_{\infty}^2)|I|,\text{ for
all intervals }I.
\end{equation}

If one is willing to apply the $T(1)$-Theorem, the proof is
complete. Otherwise, there is a simple way of avoiding $T(1)$ and reducing the estimate to the $H^1-BMO$ duality (\cite{MV}).

Before confronting the Cauchy Integral on a Lipschitz graph let us pause to discuss the relation between the Cauchy kernel and Menger
curvature, which is the main contribution of Melnikov to the subject. He was trying to find an analog for analytic capacity of
the transfinite diameter description of logarithmic capacity. In dealing with the many finite sums that 
appeared he found the following identity.
\begin{lemma*}
Given three distinct points $z_1, z_2$ and $z_3$ in the plane, one has
\begin{equation}\label{mel} 
\sum_{\sigma}\frac{1}{(z_{\sigma(2)}-z_{\sigma(1)})
(\overline{z_{\sigma(3)}-z_{\sigma(1)}})} = R(z_1,z_2,z_3)^{-2},
\end{equation}
where the sum is over all permutations of $\{1,2,3\}$ and $R(z_1,z_2,z_3)$ is the radius of the circle through $z_1, z_2$
and $z_3.$
\end{lemma*}

I tried to find an explanation to \eqref{mel} and I failed. There is no explanation: it is only a straightforward computation. 
But I learnt a way of interpreting \eqref{mel}.  If one
considers the Cauchy Integral $C$ (in the $\operatorname{p.v.}$ sense) with respect to the discrete measure 
$\delta_{z_1}+\delta_{z_2}+\delta_{z_3}$ and one computes $\|C(1)\|^2$ then one gets $\sum_{i<j}|z_i-z_j|^{-2}$ plus the left hand
side of \eqref{mel} (see \cite{V4}). This suggest that \eqref{mel} should be linked to the $T(1)$-Theorem. The most surprising feature 
of \eqref{mel} is the positivity of the left hand side. Menger had defined many years before the curvature associated with three points as the inverse of the radius of the circle through them (see \cite{V4} and \cite{K}).

Let $A$ be a Lipschitz function on $\R$ and parametrize its graph by $\gamma(x)=x+i A(x), \; x \in \R.$ The truncated
Cauchy Integral is
\[
C_{\epsilon}f(x)=\int_{|y-x|>\epsilon}\frac{f(y)\,dy}{\gamma(y)-\gamma(x))},\quad
x\in\mathbb{R},
\]
Fix an interval $I \subset \R.$ Then, by standard estimates,
\begin{equation*}
\begin{split}
\int_I|C_{\epsilon}(\chi_I)(x)|^2 \, dx&=\int_I
C_{\epsilon}(\chi_I)(x)\overline{C_{\epsilon}(\chi_I)(x)}\,dx\\*[7pt]
&=\int_I\int_{I_{\epsilon}(x)}\int_{I_{\epsilon}(x)}\frac{dy\,dz\,dx}
{(\gamma(y)-\gamma(x))(\overline{\gamma(z)-\gamma(x)})}\\*[7pt]
&=\iiint_{S_{\epsilon}}\frac{dx\,dy\,dz}{(\gamma(y)-\gamma(x))
(\overline{\gamma(z)-\gamma(x)})}+O(|I|),
\end{split}
\end{equation*}
where $I_{\epsilon}(x)$ and $S_{\epsilon}$ are as before.
By permutating the positions of the
three variables in the integrand of the right hand side of the preceding formula we get
$6$ different expressions for the left hand side. Hence
\begin{equation}\label{per}
\int_I
|C_{\epsilon}(\chi_I)(x)|^2\,dx=\frac{1}{6}\iiint_{S_{\epsilon}} R^{-2}(\gamma(x),\gamma(y),\gamma(z))\,
dx\,dy\,dz+O(|I|).
\end{equation}
By a well known formula in the geometry of the triangle $R^{-1}(\gamma(x),\gamma(y),\gamma(z))$ can be written as
$4$ times the quotient between the area of the triangle determined by $\gamma(x),\gamma(y)$ and $\gamma(z)$ and the product of the 
three side lengths. Estimating from below the length of each side by the length of its projection on the horizontal axis we get
\begin{equation*}\label{}
 R^{-1}(\gamma(x),\gamma(y),\gamma(z)) \le
4\left|\frac{\dfrac{A(y)-A(x)}{y-x}-\dfrac{A(z)-A(x)}{z-x}}{z-y}\right|,
\end{equation*}
which shows that the Cauchy Integral is controlled by the first commutator. Therefore it is bounded on $L^2(\R).$

It is disappointing that nothing similar to Menger curvature exists for kernels with homogeneity $-d$ with $ 1< d<n$ in $\R^n.$
Instead, if $0<d<1$ the symmetrization approach works perfectly well \cite{Pr}. There are still unsolved basic problems for homogeneities
with $d$ larger than $1,$ although many others have been solved using more sophisticated tools than symmetrization.

\section{David and Semmes : Uniform rectifiability}

The question is now : can we describe the measures $\mu$ with the property that the Cauchy Integral is bounded on $L^2(\mu)$ ?
The area measure on the unit disc is such a measure, but we want to stay in a context where the Cauchy kernel is no better that 
a genuine CZ
kernel. By homogeneity considerations this
amounts to requiring that the measure satisfies the AD condition

\begin{equation}\label{ad}
C^{-1}\, r \le \mu(B(z,r)) \le C \,r, \quad z \in \operatorname{supp}\mu , \quad 0< r < \operatorname{diam}(\operatorname{supp}\mu),
\end{equation}
where $C$  a positive constant independent of $z$ and $r.$  Clearly $\mu$ is boundedly absolutely continuous with respect to one 
dimensional
Hausdorff measure $H^1$ restricted to the support $E$ of $\mu.$ One can then forget about $\mu$ and inquire about the properties of the closed set  $E$ that
guarantee that the Cauchy Integral is 
$L^2(H^1|_E)$ bounded, provided $H^1|_E$ satisfies the AD condition. Of course one can ask similar questions in $\R^n$ replacing the Cauchy
kernel by other families of well behaved kernels. For example, one can fix an integer dimension $1\le d <n$ and consider  the family
$\mathcal{K}_d$ of odd kernels $K$
satisfying
\begin{equation}\label{ker}
 |\nabla^j K(x)| \le \frac{C}{|x|^{d+j}}, \quad x \neq 0, \quad j=0,1, \dots 
\end{equation}

It can be shown rather easily that if $\mu$ is $d$-dimensional Hausdorff measure on a $d$-dimensional Lipschitz graph
\begin{equation*}\label{}
\{(x_1, ..., x_d, A(x_1,..., x_d)): (x_1,...,x_d) \in \R^d \}
\end{equation*}
with $A: \R^d \to \R^{n-d}$ a Lipschitz function, then each operator with kernel in $\mathcal{K}_d$ is $L^2(\mu)$ bounded. As before,
boundedness means uniform boundedness  of the truncations
\begin{equation*}\label{}
T_\e(f)(x) = \int_{|y-x|>\e} K(x-y)\,f(y)\,d\mu(y),  \quad \e>0, \quad x \in \R^n, \quad f \in L^1(\mu).
\end{equation*}
This follows from the CMcM Theorem and the method of rotations if the kernel is homogeneous of degree $-d.$ An additional argument 
is required to get the whole family $\mathcal{K}_d$ (see \cite{DS}). However, it was known that the AD regularity condition \eqref{ad}
is not strong enough to give boundedness, even for $d=1.$ The example is the famous Garnett-Ivanov set endowed with one dimensional
Hausdorff measure. Take the unit square $Q_0=[0,1]\times [0,1].$ Divide $Q_0$ into $16$ squares of side length $1/4$ and take the 
four corner squares that contain a vertex of $Q_0.$ Repeat the process inside each of these $4$ squares and then inside the squares
one obtains at each step by this procedure. One gets at the $n$-th generation $4^n$ squares $Q_{n,j}$ of side length $1/4^n.$
The set is $E = \cap(\cup_{j=1}^{4^n} Q_{n,j}).$ There exists a unique probability measure $\mu$ on $E$ that assigns mass $1/4^n$ to each 
$Q_{n,j}.$ This measure is a positive multiple of one dimensional Hausdorff measure on $E.$ 

I claim that 
the Cauchy Integral is not $L^2(\mu)$ bounded. This can be easily seen via Menger curvature, which shows again its invaluable power. The reader should be aware that
the corresponding statement in $\R^3$ for a self-similar two dimensional Cantor set and the kernel $x/|x|^3$ is still true, but that 
one cannot
provide a so quick argument because there is nothing like Menger curvature for kernels of homogeneity $-2.$

Proceeding as in the previous section we obtain the analog of \eqref{per}
\begin{equation}\label{mcau}
 \int 
|C_{\epsilon}(\mu)(z)|^2\,d\mu(z)=\frac{1}{6}\iiint_{S_{\epsilon}} R^{-2}(z,w,\zeta
)\,
d\mu(z)\,d\mu(w)\,d\mu(\zeta) +O(\|\mu\|),
\end{equation}
where
$$
S_\e=\{(z,w,\zeta) \in E^3 : |z-w|>\e, |z-\zeta|>\e \;\text{and}\;|w-\zeta|>\e \}.
$$
If the Cauchy Integral were bounded on $L^2(\mu)$ the left hand side above would be not greater that a constant times $\|\mu\|$
and so, letting $\e$ go to $0$, we would get that
$$
\iiint_{E^3} R^{-2}(z,w,\zeta
)\,
d\mu(z)\,d\mu(w)\,d\mu(\zeta) < \infty.
$$
That this is false can be shown as follows. Choose inside each $Q_{n,j}$ three different squares of the next generation, say
$Q_{n+1,j_1}, Q_{n+1,j_2}$ and $Q_{n+1,j_3}$ and set
$$
T_{n,j}= Q_{n+1,j_1}\times Q_{n+1,j_2} \times Q_{n+1,j_3}.
$$
Since the $T_{n,j}$ are disjoint
\begin{equation}
\begin{split}
 & \iiint_{E^3} R^{-2}(z,w,\zeta
)\,
d\mu(z)\,d\mu(w)\,d\mu(\zeta)\geq \\*[7pt] 
&  \sum_{n,j} \iiint_{T_{n,j}} R^{-2}(z,w,\zeta)\,d\mu(z)\,d\mu(w)\,d\mu(\zeta)\geq \\*[7pt] 
&  c\, \sum_{n=1}^\infty \sum_{j=1}^{4^n} 4^{2n} \frac{1}{4^{3n}} = \infty. 
\end{split}
\end{equation}
The reason why the preceding argument worked is that the set $E$ is very sparse.

David and Semmes found, at first independently, examples of $d$-dimensional surfaces such that kernels in $\mathcal{K}_d$ 
are $L^2(\mu)$ bounded for $\mu= H^d,$  provided \eqref{ad} holds. In each case the boundedness could be proved from the ``Big pieces of Lipschitz graphs''
property that we described in section $4$ 
in the context of curves. A novelty emerged: Semmes introduced a new device to prove $L^2$ boundedness, namely what he called
a ``corona decomposition''.  The terminology comes from Carleson's original argument in the proof of the Corona Theorem.
To say a few words about that we start by recalling that if one has a doubling measure on a closed set $E$ then one can define families
of subsets of $E$, called dyadic squares, which enjoy the familiar properties of the dyadic net in $\R^n.$ A 
corona decomposition is a partition of the set of dyadic squares in a set of families of dyadic squares, called trees, on which your 
set $E$ can be approximated fairly well by a Lipschitz graph
associated with the tree, plus a family of bad squares where such an approximation does not happen. The family  $\mathcal{B}
$ of bad squares is small in the sense that it satisfies the Carleson condition
\begin{equation*}
 \sum_{\mathcal{B} \ni Q\subset Q_0} \mu(Q) \le C\,\mu(Q_0), \; \text{for each dyadic square}\; Q_0.
\end{equation*}
Also the squares of each tree $T$ are descendants of a top dyadic square $Q(T)$ and the family of trees is small, again in a sense
made precise by a Carleson condition:
\begin{equation*}
 \sum_{Q(T) \subset Q_0} \mu(Q(T)) \le C\,\mu(Q_0), \; \text{for each dyadic square}\; Q_0.
\end{equation*}
There is a mechanism which transports $L^2$ boundedness from the Lipschitz graphs to the measure $\mu$ on $E.$ This is, in a way, 
a sophisticated alternative of David's good lambda method to conclude from the ``Big pieces of Lipschitz graphs'' property.
It turns out that all operators with kernels in $\mathcal{K}_d$  are $L^2(\mu)$ bounded
if and only if  you have a corona decomposition. Sets having corona decompositions are called uniformly rectifiable and 
can be described by a huge number of apparently unrelated conditions, some very geometric (the reader is invited to consult \cite{DS}).
One that explains the terminology is this:
there is a large number $M$ and a small number $\theta$ such that, given a ball $B$ centered on $E,$ there exists a subset $A$ of $\R^d$
and a bilipschitz mapping
$f: A \to \R^n$ with bilipschitz constant $\le M$ with the property that $\mu( B\cap f(A))\ge \theta \mu(B).$ This condition implies
rectifiability, and is uniform at all scales and locations. David and Semmes proved that this condition follows from ``Big pieces
of Lipschitz graphs'' and still implies boundednes of all operators with kernels in $\mathcal{K}_d.$ It is strictly weaker than
``Big pieces of Lipschitz graphs'' and, in fact,  it is equivalent to having big pieces of sets which have big pieces of Lipschitz graphs.

There is a characterization of uniform rectifiability in terms of the Jones' numbers, extremely elegant. The beta number of order $1$ associated with the point $x$ and the scale $t$ is
\begin{equation*}
 \beta_1(x,t) = \inf_P \frac{1}{t^d} \int_{E\cap B(x,t)} \frac{\operatorname{dist}(y,P)}{t} \,d\mu(y),
\end{equation*}
where the infimum is taken over all affine $d$-dimensional planes  $P$  in $\R^n.$ This number measures the minimal deviation of 
the set $E$
from a $d$-dimensional plane around $x$ at the scale $t.$ Then the following Carleson condition
is equivalent to uniform rectifiability:
\begin{equation*}
\int_0^r \int_{B(x,r)}  \beta_1(y,t)^2 \, d\mu(y) \frac{dt}{t} \le C\,r^d,  \quad x \in E, \quad r>0. 
\end{equation*}
In fact, the $L^\infty$ version of the beta numbers
\begin{equation}\label{beta}
 \beta(x,t) = \inf_P  \,\frac{1}{t}\,\sup_{y\in E\cap B(x,t)}\,\operatorname{dist}(y,P) ,
\end{equation}
was introduced by Jones to give another proof of the CMcM Theorem \cite{J1}.

A difficult problem, presently known as the ``David-Semmes problem'', was left open in \cite{DS}: is the $L^2$ boundedness of Cauchy kernel ($n=2$ and $d=1$) or of 
the vector valued Riesz kernel $x/|x|^{d+1}$ for $1 \le d<n$ in $\R^n$ sufficient to guarantee uniform rectifiability ?
A first positive answer was given for $d=1$ and all $n$ in \cite{MMV} using Menger curvature. Almost $20$ years later Nazarov, Tolsa and 
Volberg solved the problem for $d=n-1$ in a formidable paper \cite{NToV1}, which  surmounts enormous difficulties and ends up with an application of the maximum principle. The problem is still 
open for integer dimensions in the range $1<d<n-1,$ because nobody has found a way of getting  around the maximum principle.

Uniform rectifiability is an extremely useful, natural, intrinsic notion that brings together the contributions of Carleson on the Corona 
Theorem and those of the CZ school via Littlewood-Paley theory. The corona decomposition makes a bridge
 between Geometry and Analysis.
At the core, again, is the magic key \eqref{cau}.

\section{Back to analytic capacity: Vitushkin's conjecture}

Vitushkin conjectured in his famous 1967 paper we have mentioned before that, for a compact set,  the property of having zero analytic 
capacity was equivalent to being Besicovitch irregular (or purely unrectifiable), that is, having projections of zero length in almost 
all directions (almost all measured with respect to arc length on the unit circle). 

In this form Vitushkin's conjecture was disproved 
by Mattila by means of a funny argument which did not allow one to decide which direction was false \cite{Ma1}. He proved that 
pure unrectifiability is not conformally invariant while $H^\infty$ removability clearly is.
Later on clever constructions
showed that a set may have positive analytic capacity and still project in zero length in almost all directions \cite{JM}, \cite{JMu}. But these sets
did not have finite length (one dimensional Hausdorff measure) and so the reformulation of Vitushkin's conjecture in which one assumes
a priori that the set has finite length remained open. For sets of finite length Besicovitch had proven that pure unrectifiability
is equivalent to intersecting 
any rectifiable curve in a set of zero length, which, by the way, explains the terminology. Now, if a compact set intersects a rectifiable curve in a set of positive length, then it also
intersects a Lipschitz graph with little Lipschitz constant in a set of positive length. Thus the 1977 theorem of Calder\'on shows
that the analytic capacity of the set is positive. 

It remained the other direction, that is, proving that if a set has positive analytic capacity
then there is a rectifiable curve that intersects the set in positive length. This looks absolutely scary. How can you construct
a rectifiable curve that intersects a given set in a set of positive length?  Peter Jones had introduced in 1987 the beta numbers 
to estimate the Cauchy Integral on Lipschitz graphs \cite{J1} . He then realized that the beta numbers could be used to define 
a geometric square function characterizing compact sets through which a rectifiable curve can pass \cite{J2}. 
Given a dyadic square $Q$ of side length $l(Q)$ set
$$
\beta_E(Q) = \inf_L  \,\frac{1}{l(Q)}\,\sup_{z \in E\cap 3Q}\,\operatorname{dist}(z,L), 
$$
where the infimum is over all straight lines $L.$  Then $2 \beta_E(Q) \,l(Q)$ is the width of the narrowest strip that contains $E\cap 3Q.$
\begin{theorem*}[\bf{P.\ Jones, 1990}]
Let $E$ be a compact subset of the plane. Then there exists a rectifiable curve $\Gamma$ such that $E \subset \Gamma$ if and only if
$$
\sum_{Q\, \text{dyadic}} \beta_E(Q)^2 \, l(Q) < \infty.
$$
\end{theorem*}

The idea for proving Vitushkin's conjecture is the following. Let $E$ be a given set of finite length. Then $E$ carries a natural
finite measure, the one dimensional Hausdorff measure restricted to $E,$  which we call $\mu.$
If one can show from the fact that $E$ has 
positive analytic capacity that the Cauchy integral is $L^2(\mu)$ bounded 
on a subset $F$ of $E$ of positive length, then one can hope that, by using Menger curvature and relating it to the 
beta numbers, $F$ can be shown to be a subset of a rectifiable curve. This was proved in 1996 in \cite{MMV} under the additional assumption
that the set $E$ is AD regular. That this scheme was plausible was suggested by a previous
result of Christ \cite{Ch}, who proved the boundedness of the Cauchy Integral on a subset of positive length under 
the assumption that the length measure on $E$ is doubling. Christ could not go farther because at that time 
Menger curvature had not yet entered the scene. The reader may consult \cite{Pa} about the Cauchy Integral, beta numbers and rectifiability.

At this point the question was: how could one get rid of the AD regularity
assumption, if $CZ$ theory had only been developed on spaces of homogeneous type? The reader should transport herself to 1996. X.Tolsa was a graduate student at Barcelona and not a word had been said on non-doubling
$CZ$ theory.  I strongly believed that it would have been
extremely useful for the understanding of analytic capacity to decide whether or not the $T(1)$-Theorem should hold for a not necessarily
doubling measure satisfying the linear growth condition $\mu(B(a,r)) \le C \, r, \; a \in \C, \; 0< r$ (which is a necessary condition
for $L^2(\mu)$ boundedness if $\mu$ has no atoms). Tolsa succeeded in proving the
non-doubling $T(1)$-Theorem for the Cauchy Singular Integral, which became part of his thesis \cite{T1}, and simultaneously Nazarov, Treil and 
Volberg, who had heard about the problem from another source, proved the same result for much more general kernels \cite{NTV1}.

After that the situation was  the following : the scheme for the proof of Vitushkin's conjecture was set up and non-doubling $CZ$ theory was rapidly developing, often overcoming  frightening difficulties posed
by non-homogeneity. Then combining two formidable papers \cite{DM} and \cite{D2} Vitushskin's conjecture was finally proved. One constructed
families of special sets
replicating the grid of dyadic squares in $\R^n$  only under the linear growth hypothesis, and 
a special ``ad hoc'' non-homogeneous $T(b)$-Theorem was shown to hold, which provided a piece $F$ of positive length on which 
the Cauchy Integral was $L^2(\mu)$ bounded. Then the Menger curvature of the length measure on $F$ 
turned out to be finite. This implied rectifiability by a theorem proved by David in the plane and extended to several dimensions by L\'eger \cite{Le}. In this last step a corona
type decomposition played a central role, as well as the beta numbers, which were finally responsible for landing on rectifiability. 
See the survey \cite{D3} for an agile presentation of the main ideas.

The proof of Vitushkin's conjecture was an impressive tour de force that encompassed many remarkable contributions. 
The reader should note that most of the tools used were fabricated in the environment of the CMcM Theorem : $T(b),$ the beta
numbers, Menger curvature, corona decompositions. 

It is worthwhile at this point to make a historical remark. The sequence of events concerning Menger curvature was as follows. First
Melnikov discovered the relation between the Cauchy kernel and Menger curvature \eqref{mel}. Then the author found the proof of the CMcM Theorem  discussed in section 6 in \cite{MV} (see also \cite{V2}). Then with Mattila we solved the one dimensional David-Semmes problem and the homogeneous 
Vitushkin's conjecture in \cite{MMV} and after that main result of  \cite{Me} was proved. In \cite{Me} the other two papers,
which already existed in preprint form, were not mentioned. 

Another approach to Vitushkin's conjecture was devised by Nazarov, Treil and Volberg, just after the David-Mattila proof, but 
it has
never been published in article form in a journal, for reasons unknown to the author. A draft was circulated at the time and its
contents have been partially 
presented in the books \cite{Du}, \cite{Vo} and \cite{TB}. A recent version has been uploaded to the ArXiv by the authors \cite{NTV2}. 
This magnificent paper has 
been very influential. For instance, it lead Tolsa to a proof of the semiaddditivity of analytic capacity and to a solution of 
Painlev\'e's problem. It has also been used in other deep problems concerning rectifiability.

\section{The semiadditivity of analytic capacity and the \\solution of Painlev\'e's problem}
The NTV(Nazarov-Treil-Volberg) proof of Vitushkin's conjecture is based on a special $T(b)$-Theorem. One has a compact set $E$ of finite 
length (one dimensional Hausdorff measure)
and positive analytic capacity. Since $E$ is any set with these properties, the length measure $\mu$ on $E$ is not necessarily doubling. The fact that the set has positive analytic capacity yields
easily a complex $\mu$-measurable
bounded function $b$ on $E$ with $\int b \, d\mu >0$ and bounded Cauchy Integral. Of course we would like to apply the $T(b)$-Theorem,
but at least two enemies are barring our way. 
The first is that for 
almost all $z\in E$ one has $\mu(D(z,r)) \le C \,r, \; r < r_0,$  with $C$ and $r_0$ depending on $z,$ where $D(z,r)$ is the disc centered at $z$ of
radius $r.$  Then, given $M>0$ there might be many non-Ahlfors discs
satisfying $\mu(D(z,r)) > M \,r.$ The second enemy is that we do not have the para-accretivity condition $|\int_D b \,d\mu | > c\, \mu(D)$
for all discs $D.$ By essentially taking the union of the non-Ahlfors discs and the non-para-accretive discs one forms an open set $U$ 
on which all bad things happen. Then one ingeniously modifies the Cauchy kernel on $U$ and through a sophisticated lengthy sequence of 
clever moves, sometimes hard,
one ends up proving that the modified kernel provides an operator which is bounded on $L^2(\mu)$. Since the Cauchy kernel has remained 
untouched on $F=E \setminus U$ the Cauchy Integral is a bounded operator on 
$L^2(\mu|_F).$  But one has taken care to ensure that $\mu(U)$ is small. Hence $F$ has positive length and you are done, after 
remarking that the Menger curvature of $\mu$ restricted to $F$ is finite and applying
David-L\'eger.
 
One uses several pieces of the well-known tool box, like martingale differences associated with $b,$ in the spirit of \cite{CJS}, Schur's lemma, Carleson type
conditions and a beautiful probabilistic argument on random grids of standard dyadic squares (which was first devised in \cite{GJ}).

What is the relation of $T(b)$ with the semiadditivity of analytic capacity ? To have a glimpse of that we first discuss a 
variant of analytic capacity, called positive analytic capacity, defined as
\begin{equation*}\label{}
 \gamma^+(E) = \sup \mu(E),
\end{equation*}
where the supremum is taken over all positive measures supported on $E$ such that $\frac{1}{z}  \star \mu$ is a function in
$L^\infty(\C)$ with norm less than or equal to $1.$ Clearly
\begin{equation*}\label{}
 \gamma^+(E) \le  \gamma(E).
\end{equation*}
To better understand positive analytic capacity set
\begin{equation}\label{gamop}
 \gamma_{\operatorname{op}}(E) = \sup \mu(E),
\end{equation}
where the supremum is over all positive measures supported on $E$ satisfying
\begin{equation}\label{lgrowth}
\mu(D(z,r)) \le \,r, \quad z \in \C, \quad r>0,
\end{equation}
and such that the operator norm of the Cauchy Integral as a bounded operator from $L^2(\mu)$ into itself is less than or equal to $1.$

It is not obvious, but not difficult, to show that the non-doubling $T(1)$-Theorem for the Cauchy Singular Integral yields
\begin{equation*}\label{}
\gamma^+(E) \le C\, \gamma_{\operatorname{op}}(E).
\end{equation*}
On the other hand, the fact that standard $CZ$ theory works also in a non-doubling context plus the Davie-Oksendal dualization method of a 
weak-type inequality gives easily the reverse inequality. This is important for semiadditivity. Indeed, take two compact
sets $E_1$ and $E_2.$  It is clear 
from the definition that if $\mu$ is an admissible measure for the supremum defining $\gamma_{\operatorname{op}}$ relatively 
to $E_1 \cup E_2$
then $\mu_j = \mu|_{E_j}$ are  admissible measures for $E_j$. Hence
\begin{equation*}\label{}
\gamma_{\operatorname{op}}(E_1 \cup E_2) \le \gamma_{\operatorname{op}}(E_1)+\gamma_{\operatorname{op}}(E_2).
\end{equation*}
Being comparable to a subadditive set function one concludes that $ \gamma^+$ is semiadditive. Note that this fact is not at all obvious
from the definition, the reason being that if one restricts a positive measure $\mu$ such that  $\frac{1}{z}  \star \mu$ is a function in
$L^\infty(\C)$ to a general compact $F \subset \operatorname{spt}\mu$ is not true any more that $\frac{1}{z}  \star \mu|_F$ is in $L^\infty(\C).$ 

To prove semiadditivity of analytic capacity it remains to show that $\gamma$ is comparable to $\gamma^+$ or, in other words, that 
there exist a constant $C_0$ such that
\begin{equation}\label{gamsem}
\gamma(E) \le C_0 \, \gamma^+(E), \quad \text{for all compact}\; E \subset \C.
\end{equation}
This is what Tolsa proved in \cite{T1} using the NTV $T(b)$-Theorem. The proof features several subtle points, some rather involved, and  the remarkable technical power of the author leaves a deep impression on the reader. 

As a corollary one obtains a solution of Painlev\'e's problem. If $E$ has positive analytic capacity then there exists a positive
measure $\mu$ supported on $E$ with bounded Cauchy potential. Then $\mu$ has linear growth \eqref{lgrowth} and by \eqref{mcau} 
finite Menger curvature, 
namely,
\begin{equation}\label{pain}
\iiint R^{-2}(z,w,\zeta)\,d\mu(z)\,d\mu(w)\,d\mu(\zeta) < \infty.
\end{equation}
The converse is also true \cite{Me} and is now a straightforward consequence of the non-doubling $T(1)$-Theorem for the Cauchy Singular Integral \cite{V3}. If $\mu$ is a non-zero positive measure supported on $E$ satisfying the linear growth condition
\eqref{lgrowth} and \eqref{pain} then it follows readily from the $T(1)$-Theorem that the Cauchy Integral is $L^2(\mu|F)$ bounded for some $F$ with positive
$\mu$ measure. By CZ theory and Davie-Oksendal the analytic capacity of $E$ is positive. 

Condition \eqref{pain} has certainly a geometric flavour and consequently, with some generosity, one could say that is a solution to
Painlev\'e's problem. I asked Tolsa whether \eqref{pain} was really a geometric condition in the sense that it was  a bilipschitz invariant and his  (positive) answer was the paper \cite{T4}. In fact, whether analytic capacity is, modulo constants, a bilipschitz invariant is a question raised by O'Farrell (see, for instance, \cite{DOF}, where the question was whether analytic capacity is a linear invariant, modulo constants).

It is a challenge to describe the main ideas of the proof of \eqref{gamsem}, because one has to let technical issues pervade the exposition. 
Nevertheless, I cannot refrain from mentioning one of these ideas in a particular case \cite{MTV}. The uninterested reader may proceed happily to the next section.

In proving \eqref{gamsem} one can assume that $E$ is a union of finitely many squares. Let us consider, as a significant example,  the union $E_N$ of the $4^N$ 
squares
$Q^N_j, \; 1 \le j\le 4^N,$ of generation $N$ in the construction of the corner quarters Cantor set described in section 7. The side length 
of $Q^N_j$ is $4^{-N}.$ It is known that $\gamma^+(E_N) \approx N^{-1/2}$ (\cite[p. 125]{TB}). Take a function $f,$ analytic off
$E_N,$ with $|f(z)|\le 1, \; z \notin E_N$ and $f'(\infty) = \gamma(E_N).$ Then
\begin{equation*}\label{}
f(z)= \frac{1}{2\pi i} \int_{\partial E_N} \frac{f(\zeta)}{z-\zeta}\,d\zeta, \quad z \notin E_N,
\end{equation*}
by Cauchy's integral formula. Thus $f=C(\nu),$ where $\nu$ is the complex measure $\frac{f(\zeta)\,d\zeta}{2\pi i}.$ Our goal is 
to find a positive measure $\mu$ supported on $E_N,$ with mass $\gamma(E_N),$ such that the norm of the Cauchy Singular Integral, as a bounded operator on $L^2(\mu),$ 
is bounded by an absolute constant. Such absolute bound should follow from an application of some kind of $T(b)$-Theorem. It is here
very tempting to try Christ's local $T(b)$-Theorem, which states, in the present context,  that $L^2(\mu)$ boundedness follows if for each square 
$Q^n_j, \; 1\le j \le 4^ n, \; 1 \le n \le N,$ one finds a function $b^n_j,$ \; $|b^n_j|\le 1,$ with $|C(b_j\mu)|\le 1$ off $E_N,$ satisfying
the local para-accretivity condition $|\int_{Q^n_j} b^n_j \,d\mu | \ge c\,\mu(Q^n_j).$
In other words, instead of looking for one global para-accretive bounded $b,$ with bounded $C(b\,d\mu),$ one has the flexibility of 
making a local construction (with absolute constants). One starts with the naive choice
\begin{equation*}\label{}
d\mu = \frac{1}{4}\gamma(E_N) ds_{\partial E_N}.
\end{equation*}
Given a square $Q^n_k, \; 1 \le n \le N,\; 1\le k\le 4^N,$ one sets, as a first try for $b^n_k,$
\begin{equation}\label{betan}
\beta^n_k = \sum_{Q^N_j \subset Q^n_k} \frac{1}{2\pi i} f(\zeta)\,d\zeta_{\partial Q^N_j}.
\end{equation}
It is a simple result, part of Vitushkin's localization methods, that
\begin{equation}\label{cbetan}
|C(\beta^n_k)(z)| \le A, \quad z \notin Q^n_k \cap E_N,
\end{equation}
for some absolute constant $A.$ But
\begin{equation*}\label{}
|\beta^n_k|= \sum_{Q^N_j \subset Q^n_k} \frac{1}{2\pi} |f(\zeta)| \,ds_{\partial Q^N_j}
= \frac{4}{\gamma(E_N)} \sum_{Q^N_j \subset Q^n_k} \frac{1}{2\pi} | f(\zeta) | \, d\mu(\zeta),
\end{equation*}
which blows up with $N$ because we know that $\gamma(E_N)\rightarrow 0$ as $N \to \infty.$
One has to modify the naive approach and the right change is a little subtle manouver.
 
Define a new $\mu$ by 
$$\mu=\gamma(E_N) \frac{ds_{\partial E_M}}{\operatorname{length}(E_M)}$$
for some $M$ between $1$ and $N,$  which has to be chosen.
Note that our new $\mu$ has support inside $E_M,$ which is larger than $E_N.$ Then, in the event that we could show that
the Cauchy Singular Integral, as a bounded operator from $L^2(\mu)$ into itself, has norm bounded by an absolute constant we would get
$\gamma(E_N)=\mu(E_M) \le C_0\, \gamma^+(E_M).$  Hence we need $\gamma^+(E_M) \le C_0 \, \gamma^+(E_N).$ Since 
$\gamma^+(E_N) \approx N^{-1/2} $ it is enough to take $M \ge N/2$. Indeed the right choice, as we will check now,  is $M=N/2$ (assume $N$ even). 

Take a new $\beta^n_k, $  for $ 1 \le n \le M$  and $1 \le k \le 4^n,$ as in \eqref{betan} with $N$ replaced by $M.$  We still have 
the favourable bound \eqref{cbetan} but again the maximum of $|\beta^n_k |$ blows up with $N.$ To modify $\beta^n_k$ we first remark
that there are simple ways of constructing smooth functions $\varphi^n_j$ supported on $Q^n_j$
such that $ 0\le \varphi^n_j \le 1,$ $|C(\varphi^n_j\,ds_{\partial Q^n_j})|\le A$ and $\int \varphi^n_j \,ds_{\partial Q^n_j} \ge \frac{1}{4^n}.$
One just dilates and translates a model function on $Q_0=[0,1]\times [0,1] .$   It is useful in modifying $\beta^n_k$ to preserve the
integral on each building piece $Q^M_j,$ which is $\nu(Q^M_j).$ Therefore we set
\begin{equation*}\label{}
b^n_k = \sum_{Q^M_j \subset Q^n_k} \nu(Q^M_j) \frac{\varphi^M_j}{\int \varphi^M_j \,d\mu}, \quad 1\le n\le M, \quad 1\le k \le 4^n.
\end{equation*}
Clearly, on the one hand, we have  $\int \varphi^M_j \,d\mu \ge \gamma(E_N) \,\frac{1}{4^M} = \mu(Q^M_j).$  On the other hand,
$C(\chi_{Q^M_j} \nu) $ is analytic off $Q^M_j \cap E_N$ and
\begin{equation*}\label{}
|C(\chi_{Q^M_j} \nu )(z)| \le A, \quad  z \notin Q^M_j \cap E_N,
\end{equation*}
by Vitushkin's localization technique. Thus, by definition of analytic capacity,  
$$|\nu(Q^M_j)| \le A 	\, 	\gamma(Q^M_j \cap E_N).$$ 
%\vspace{2cm}
Now $Q^M_j \cap E_N$ is the result of applying $N-M=N/2$ steps of the construction
of the corner quarters Cantor set, starting from $Q^M_j.$  Hence, by the homogeneity of analytic capacity, $\gamma(Q^M_j \cap E_N)=
1/4^{N/2} \,\gamma(E_{N/2}).$  At this point one should note that we do not know the precise relation between $\gamma(E_{N})$ and 
$\gamma(E_{N/2}),$ besides the obvious fact that $\gamma(E_{N}) \le \gamma(E_{N/2}).$ Assume for a moment that $\gamma(E_{N/2}) 	\le C_1 \, \gamma(E_{N}),$ for some large absolute constant $C_1$ to be chosen later. Then $\nu(Q^M_j) 
\le A\,\gamma(Q^M_j \cap E_N) \le A C_1 1/4^M \,\gamma(E_N) = A C_1 \mu(Q^M_j),$ which tells us that $b^n_k $ is bounded by an absolute constant. Note that we have overcome the main difficulty which appeared  in the naive approach we started with.

A standard argument, based on the quadratic decay of each term, allows to control the supremum norm of the difference
$C(\beta^n_k) - C(b^n_k)$ and thus $|C(b^n_k)(z)| \le A, \; z \notin Q^n_k \cap \partial E_M.$  

It remains to check the local para-accretivity condition.  Since $\sum_{j=1}^{4^M} \nu(Q^M_j) = \gamma(E_n),$ for at least one index $k$ one has $|\nu(Q^M_k)| \ge \gamma(E_N)/4^M=\mu(Q^M_j).$
For this special $k$ the construction of the function $b^n_k$ required by Christ's local $T(b)$-Theorem is complete. One may wonder what one should do for $Q^n_j$ with $j \ne k.$ The answer is simple : transport the function $b^n_k$ by translation. The conclusion is that in this case \eqref{gamsem} follows from an application of the local $T(b)$-Theorem.

It remains to deal with the case $\gamma(E_N) \le C_1^{-1} 	\, \gamma(E_{N/2}).$ Proceeding by induction
\begin{equation*}\label{}
\gamma(E_N) \le C_1^{-1} 	\, \gamma(E_{N/2}) \le C_1^{-1}\,C_0  \gamma^+(E_{N/2}) \le A \,C_1^{-1}\,C_0  \gamma^+(E_{N}),
\end{equation*}
and then it suffices to choose $C_1=A.$

In the general case of the proof of \eqref{gamsem} one also chooses a union of squares containing $E,$ close enough to $E$
so that $\gamma^+$ is not altered too much, but far enough so that appropriate estimates can be carried over. Also the local $T(b)$-theorem is replaced by the NTV $T(b)$-theorem.

What we said for the proof of Vitushkin's conjecture applies here. The proof of the semiadditivity of analytic capacity and the solution to Painlev\'e's problem are extraordinary results, which were built over previous brilliant work by excellent mathematicians. The NTV $T(b)$-theorem makes the connection to the magic key \eqref{caul2}. In the next section we will review applications of $L^2$ estimates to rectifiability.

\section{The Cauchy Singular Integral and rectifiability}
It is a standard fact in classical  (first generation : convolution smooth homogeneous) $CZ$ theory that after proving the $(1,1)$ weak-type estimate for the maximal singular integral of a finite Radon measure, one can readily show the existence of the principal values of the singular integral of a finite Radon measure \text{a.e.}\ with respect to Lebesgue measure. This is still true for the Cauchy Singular Integral of a finite Radon measure $\mu$ supported on a rectifiable curve, by the Calder\'on theorem on Lipschitz graphs with small constant, as we discussed in section 2 for functions integrable with respect to $ds.$ In other words, one has the $ds$ \text{a.e.}\ existence of 
\begin{equation}\label{vpmu}
\operatorname{p.v.} C(\mu)(z) = \lim_{\e \to 0} C_\e(\mu)(z),
\end{equation}
where
\begin{equation*}
C_\e(\mu)(z)= \int_{|\zeta-z |>\e} \frac{d\mu(\zeta)}{\zeta-z}.
\end{equation*}

Although the result above is a consequence of $L^2(ds)$ estimates on Lipschitz graphs (with small constant) and thus a great achievement, it is also, in a sense, intuitively plausible. At a point $z_0$ with a tangent the curve looks symmetric with respect to $z_0$ and therefore the oddness of the Cauchy kernel can play a role in cancelling terms in opposite sides. In the other direction, if the set has finite one dimensional Hausdorff measure and is irregular in the sense of Besicovitch (like the corner quarters Cantor set of section 7), then one should think of the set as dispersed irregularly around the point and cancellation may not occur, thus making problematic the existence of the principal value.

Pertti Mattila discovered, just before the David-Semmes theory of uniform rectifiability appeared, that the intuition just described is correct  \cite{Ma2}. He proved that if a finite positive Radon measure $\mu$ has principal values $\mu$ \text{a.e.}\ and the lower density 
\begin{equation*}
\liminf_{r \to 0} \frac{\mu D(z,r)}{r} 
\end{equation*}
is positive $\mu$ \text{a.e.}\ then  $\mu$ is rectifiable in the sense of Federer, namely,  lives in a countable union of rectifiable curves. In particular,
if $\mu$ is the length measure (one dimensional Hausdorff measure) on a set $E$ of finite length, then $E$ is rectifiable if and only if 
the principal values \eqref{vpmu} exist $\mu$ \text{a.e.} and the lower density is positive.

Some hypothesis on the positivity of density is necessary: if $\mu$ is the $2$-dimensional Lebesgue measure on a ball, the principal values do exist \text{a.e.} but the measure is not rectifiable. Mattila's proof goes by tangent measures, a notion introduced by Preiss to prove rectifiability from the existence of density (in higher dimensions).  The problem of replacing the positivity of the lower density by that of the upper density remained open. It was solved in the positive in \cite{T5} using the NTV $T(b)$-theorem. Indeed, existence of principal values was replaced by $\mu$ \text{a.e.}\  finiteness of the maximal Cauchy Integral 
\begin{equation*}
C^*(\mu)(z) = \sup_{\e >0} |C_\e(\mu)(z)|
\end{equation*}
and the positivity of the lower density was replaced by the $\mu$ \text{a.e.}\ positivity of the upper density. Under these assumptions $\mu$ can be written as the sum of a discrete measure (i.e., with countable support) plus a measure
absolutely continuous with respect to length on a countable union of rectifiable curves. The $L^2$ estimates provided by NTV yield, by the usual path (finite Menger curvature and the David-L\'eger rectifiability criterion), the rectifiable curves. In particular, there are no non-zero continuous singular measures $\mu$ on the line with principal values $\mu$ \text{a.e.}\ 
Another beautiful consequence concerns the case of sets $E$ in the plane with finite length, in which no additional hypothesis on density
is required. It turns out that $E$ is rectifiable if and only if the length measure $\mu$ has principal values $\mu$ \text{a.e.}\  or if and only $C^*(\mu)(z)< \infty$ $\mu$ \text{a.e.}\  This follows from the fact that the upper density is positive and finite $\mu$ \text{a.e.}\  and so the discrete part vanishes.

The main result of  \cite{T5} is a general structural theorem for a positive finite Radon measure $\mu$ with the property that the maximal Cauchy Integral is finite $\mu$ \text{a.e.}\  Such a measure can be written as the sum of three measures $\mu=\mu_d +\mu_r+\mu_0,$ where $\mu_d$ is discrete, $\mu_r$  lives in a countable union of rectifiable curves and is absolutely continuous with respect to length and $\mu_0$  is a sum of a sequence of measures with $0$ linear density and finite Menger curvature. The measure $\mu_0$  can also be described, without mentioning Menger curvature, as the sum of a sequence of measures of $0$ linear density on which the Cauchy Singular Integral is a bounded operator on $L^2$ of the measure. As far as I know, there is no extension to higher dimensions, very likely because the proof depends on symmetrisation arguments of the kind that lead to Menger curvature and those work only for kernels with homogeneity $-\alpha$ with $0<\alpha \le 1.$

In higher dimensions there is an analog of Mattila's result, due to Mattila and Preiss \cite{MP}, in which the Cauchy kernel is replaced by the vector valued Riesz kernel $x/|x|^{n}$ of homogeneity $-(n-1)$ (modulo constants, it is the gradient of the fundamental solution of the Laplacian). Again there is a hypothesis of positivity of the  lower density in dimension $n-1.$  Tolsa extended that result omitting the lower density requirement. The statement is as follows. 

Let $\mu$ be a positive finite Radon measure and let $E$ be the set where the upper density in dimension $n-1$ is positive and finite
and, in addition, the principal value of the Riesz transform exists. In other words, $E$ is the set of $x \in \Rn$ such that
\begin{equation*}
0< \limsup_{r \to 0} \frac{\mu B(x,r)}{r^{n-1}} < \infty \quad \text{and there exists} \quad \lim_{\e \to 0} \int_{|y-x|>\e} \frac{y-x}{|y-x|^n}\,d\mu(y).
\end{equation*}
Then $E$ is rectifiable, that is, is contained in a set of the form $Z \cup (\cup_j S_j)$ where $Z$ has vanishing $(n-1)$-dimensional Hausdorff
measure and each $S_j$ is a $C^1$ hypersurface of dimension $n-1.$ 

In particular, if $\mu=H^{n-1}$ is $(n-1)$-dimensional Hausdorff measure on a set $E$ with $H^{n-1}(E) < \infty,$ then $E$ is rectifiable if and only if
the principal value 
\begin{equation*}
 \lim_{\e \to 0} \int_{|y-x|>\e} \frac{y-x}{|y-x|^n} \,dH^{n-1}(y)
\end{equation*}
exists $H^{n-1}$ \text{a.e.}\
I am mentioning this higher dimensional extension of the planar result, which is extremely valuable in itself, because of the method of proof. In the plane one uses Menger curvature
and the David-L\'eger rectifiability criterion. In higher dimensions Menger curvature is replaced by sharp estimates of the norm of the Singular Riesz transform as a bounded operator on the $L^2$ space of the restriction of $H^{n-1}$ to an $(n-1)$-dimensional Lipschitz graph. Then these estimates are combined with Corona type decomposition in the spirit of David-L\'eger.  So we see here in action again the power
of $L^2$ estimates for basic singular integral operators on Lipschitz graphs. 

One may wonder what happens with the condition of finiteness of the maximal Riesz transform:
\begin{equation*}
 %R^*(\mu)(x) 
 \sup_{\e >0} \left|\int_{|y-x|>\e} \frac{y-x}{|y-x|^n} \,d\mu(y)\right|, \quad x \in \Rn.
\end{equation*}
It can be shown \cite{NToV2} that this condition combined with the $\mu$ \text{a.e.}\ finiteness of the upper $(n-1)$-dimensional density implies rectifiability
of $\mu,$  that is, $\mu$ has rectifiable support and is absolutely continuous with respect to $H^{n-1}.$  This follows from the deep difficult main result of \cite{NToV1}: if $\mu$ is Ahlfors-David regular of dimension $n-1,$  then $L^2(\mu)$ boundedndess of the Riesz vectorial transform of homogeneity $-(n-1)$ implies uniform rectifiability.  It is worth mentioning that the preceding 
statement for Riesz transforms of integer homegeneity $-m$ with $1 < m < n-1$  seems out of reach presently.

The reader interested in the general theory of rectifiability (and uniform rectifiability) will find in the literature of the last decade a wealth of elegant criteria, most of them depending on various square functions with a geometric content. Variants of the beta numbers play a central role here. There are also interesting square functions built in terms of local approximations of a measure in terms of ``flat" measures. These approximations are measured in terms of distances in the set of Radon finite measures (distances of Wasserstein type). The papers \cite{ATT}, \cite{Da} and \cite{B} contain relevant results and an extensive list of references to many other important works, including applications to the singular sets of solutions of certain PDE.

There are two specially important recent advances involving rectifiability. The first is the solution of the $\e^2$ conjecture of Carleson
by Jaye, Tolsa and Villa \cite{JTV}. The second is the clarification role played by rectifiability in understanding under which conditions harmonic measure in $\Rn$ is absolutely continuous with respect to $H^{n-1}.$ This will be briefly reviewed in the next and last section.

\section{The Cauchy and Riesz Singular Integrals and harmonic measure}
Perron devised a general way of solving the Dirichlet problem on a domain $\Omega$ in $\Rn.$ Given a continuous function $f$ on 
$\partial \Omega,$ let $u_f$ stand for  the supremum of subharmonic functions $u$ on $\Omega$ which lie below $f$ at the boundary, namely, such that
$$
\limsup_{x \to a} u(x) \le f(a), \quad a \in \partial \Omega.
$$
Then $u_f$ is harmonic on $\Omega$ and, under mild hypotheses on $\Omega,$  takes continuously the boundary values $f.$ Hence $u_f$ solves the Dirichlet problem with boundary data $f.$ The mapping $f \rightarrow u_f$ is linear and, for each fixed $x \in \Omega,$
$f \rightarrow u_f(x)$ is linear and bounded on the space of continuous functions on the boundary. By the Riesz representation theorem this mapping is given by integrating $f$ against a positive Radon measure $\omega_x$ on $\partial \Omega:$ 
\begin{equation*}
u_f(x) = \int_{\partial \Omega} f(y)\,d\omega_x(y), \quad f \in C(\partial \Omega).
\end{equation*}
The measure $\omega_x$ is called harmonic measure with pole $x.$
By Harnack's principle, for $x, y \in \Omega$ the measures $\omega_x$ and $\omega_y$ are mutually absolutely continuous and so for many questions concerning harmonic measure the pole can be kept in the background,  so that one uses the notation $\omega$ for $\omega_x$.

More than one century ago the Riesz brothers proved that if $\Omega$ is a Jordan domain in the plane with rectifiable boundary, then $\omega$ is absolutely continuous with respect to arc length on the boundary. Since then there has been an enormous amount of work on relating geometric properties of the boundary with properties of harmonic measure (see \cite{GM} and \cite{Tr}). Bishop and Jones showed that for a planar simply connected $\Omega$ the conclusion of the F. and M. Riesz  theorem holds on a rectifiable piece of the boundary, but fails for non simply connected domains \cite{BJ}. In higher dimensions the problem becomes much more complicated. Nevertheless, the relation between $\omega$ and the $(n-1)$-dimensional 
Hausdorff measure on the boundary has been clarified recently in a seven authors paper \cite{A7}, which culminates a series of efforts by many people. The result reads as follows. 

Let $\Omega$ be a domain in $\Rn$ and $E$ a subset of $\partial \Omega$ with $H^{n-1}(E) < \infty.$  If $\omega$ is absolutely continuous with respect to $H^{n-1},$ then $\omega$ lives on a countable union of $C^1$ hypersurfaces of dimension $n-1$  (and so $\omega$ is rectifiable.)  If $H^{n-1}$ is absolutely continuous with respect to $\omega$ on $E,$ then $E$ is rectifiable.

I mention this particular result among many others which have been proved in the last decade, all extremely interesting, because there is a connection with the $L^2$ boundedness of the Singular Riesz transform (the Cauchy Singular Integral in the plane). The connection between harmonic measure and the Riesz transform is very simple, but apparently it had never been explicitly exploited before. 

Let us restrict our attention to the plane and, for the sake of convenience, assume that $\Omega$ is the complement in the Riemann
sphere of a compact set $K$ and that $\omega$ is harmonic measure with respect to the point at $\infty.$  The gradient of the logarithmic potential of $\omega$  is given by a constant times the conjugate Cauchy integral of $\omega.$ This is the connection between harmonic measure and the Cauchy Integral alluded to before. By standard properties of harmonic measure one gets the estimate, subtle but not extremely difficult to prove,
\begin{equation}\label{est}
\limsup_{\e \to 0} C_\e(\omega)(z) \le C \, \theta^*(z), \quad  z \in K, \; \omega \; \text{a.e.}
\end{equation}
where $C_\e$ is the truncated Cauchy Integral and $\theta^*(z)= \limsup_{r \to 0} \omega D(z,r)/ r$ the upper linear density.
If $\omega$ is absolutely continuous with respect to $H^1$ on $E \subset K,$  then $\theta^*$ is finite $\omega$ \text{a.e.}\ and, by
\eqref{est},  the maximal Cauchy Integral is finite $\omega$ \text{a.e.}\ on $E.$ Therefore the rectifiability criterion discussed in the previous section applies and concludes the proof (the statement about absolute continuity of $H^1$ with respect to $\omega$ follows easily). Let us remark that the rectifiability criterion we are applying depends on the NTV $T(b)$ theorem and in higher dimensions one appeals to the solution of the David-Semmes problem in codimension $1$ \cite{NToV1}.

\section{Conclusion} 
Many things have happened since 1982.  In Classical Analysis long-standing open problems, which looked inaccessible at that time, are now reasonably well understood.  New sophisticated tools and ideas have been introduced to solve them and are now used in other fields. In particular, the magic key has opened several solid doors and new landscapes have been offered to the contemplation of analysts of all tribes. The future is here, in front of us, and is bringing new challenges, new problems that hopefully will be solved, following the guidance of our predecessors.

\vspace{0.5cm}
\begin{acknowledgements} 
I would like to acknowledge the valuable assistance of various colleagues in correcting a first draft of the paper: 
J.Cuf\'i, G.David, J.B. Garnett, A.G. O'Farrell, S.Semmes and X.Tolsa.
The author is partially supported by grants
2017-SGR-395 (Generalitat de Cata\-lunya), MTM2016-75390 (Mineco) and  PID2020-112881GB-I00. 
 \end{acknowledgements}

\vspace{0.5cm}
{\small
\begin{tabular}{@{}l}
J.\ Verdera \\ Departament de Matem\`{a}tiques, Universitat Aut\`{o}noma de Barcelona\\
and\\
Centre de Recerca Matem\`atica\\
Barcelona, Catalonia\\
{\it E-mail:} {\tt jvm@mat.uab.cat}\\*[5pt]

\end{tabular}}
\end{document}